\documentclass[10pt]{article}
\usepackage[hmargin=0.9in,vmargin=0.8in]{geometry}
\usepackage[T1]{fontenc}
\usepackage{lmodern}

\usepackage{amsmath,amssymb,amsthm}
\usepackage{tikz-cd}
\usepackage{microtype}
\usepackage{needspace}
\usepackage[hidelinks]{hyperref}

\makeatletter
\g@addto@macro\normalsize{%
  \setlength{\abovedisplayskip}{6.5pt plus 2pt minus 1pt}%
  \setlength{\belowdisplayskip}{6.5pt plus 2pt minus 1pt}%
  \setlength{\abovedisplayshortskip}{3pt plus 1pt}%
  \setlength{\belowdisplayshortskip}{5pt plus 2pt minus 1pt}%
}
\g@addto@macro\th@plain{\thm@preskip=5pt plus 1pt minus 1pt\relax\thm@postskip=5pt plus 1pt minus 1pt\relax}
\g@addto@macro\th@definition{\thm@preskip=5pt plus 1pt minus 1pt\relax\thm@postskip=5pt plus 1pt minus 1pt\relax}
\g@addto@macro\th@remark{\thm@preskip=5pt plus 1pt minus 1pt\relax\thm@postskip=5pt plus 1pt minus 1pt\relax}
\renewcommand\section{\@startsection{section}{1}{\z@}%
  {-2.8ex \@plus -0.8ex \@minus -.2ex}%
  {1.3ex \@plus .2ex}{\normalfont\Large\bfseries}}
\renewcommand\subsection{\@startsection{subsection}{2}{\z@}%
  {-2.2ex \@plus -0.6ex \@minus -.2ex}%
  {.9ex \@plus .2ex}{\normalfont\large\bfseries}}
\renewcommand\paragraph{\@startsection{paragraph}{4}{\z@}%
  {1.8ex \@plus .6ex \@minus .2ex}%
  {-1em}{\normalfont\normalsize\bfseries}}
\renewcommand\@maketitle{%
  \newpage\null\vskip .7em
  \begin{center}%
    \let\footnote\thanks
    {\LARGE\@title\par}%
    \vskip 1em
    {\large\lineskip .5em
      \begin{tabular}[t]{c}\@author\end{tabular}\par}%
    \vskip .5em
    {\normalsize\@date}%
  \end{center}%
  \par\vskip .8em
}
\let\originalthebibliography\thebibliography
\let\originalendthebibliography\endthebibliography
\renewenvironment{thebibliography}[1]
  {\originalthebibliography{#1}\small
   \setlength{\itemsep}{3pt plus 1pt}\setlength{\parsep}{0pt}}
  {\originalendthebibliography}
\makeatother

\title{Non-admissible smooth irreducible representations of $\GL_2(F)$ for totally ramified extensions $F/\Qp$}
\author{Dayu Huang}
\date{November 3, 2026}

\theoremstyle{plain}
\newtheorem{theorem}{Theorem}[section]
\newtheorem{proposition}[theorem]{Proposition}
\newtheorem{lemma}[theorem]{Lemma}
\newtheorem{corollary}[theorem]{Corollary}

\theoremstyle{definition}
\newtheorem{definition}[theorem]{Definition}

\theoremstyle{remark}
\newtheorem{remark}[theorem]{Remark}

\newcommand{\Qp}{\mathbb{Q}_p}

\newcommand{\Fp}{\mathbb{F}_p}

\newcommand{\Fpbar}{\overline{\mathbb{F}}_p}
\newcommand{\GL}{\mathrm{GL}}
\newcommand{\Ind}{\mathrm{Ind}}

\newcommand{\Mat}{\mathrm{M}}

\newcommand{\Inj}{\mathrm{Inj}}
\newcommand{\Sym}{\mathrm{Sym}}
\newcommand{\Hom}{\mathrm{Hom}}
\newcommand{\soc}{\mathrm{soc}}

\newcommand{\End}{\mathrm{End}}
\newcommand{\Id}{\mathrm{Id}}

\begin{document}
\maketitle

\begin{abstract}
Let $p\geq 5$, and let $F/\Qp$ be a nontrivial finite totally ramified extension. We construct a non-admissible smooth irreducible representation of $\GL_2(F)$ over $\Fpbar$. The construction uses representations of $\GL_2(\mathcal O_F)$ on which the first principal congruence subgroup acts nontrivially, together with a deformation of a diagram indexed by $\mathbb Z$.
\end{abstract}

\section{Introduction}

Let $F$ be a finite extension of $\Qp$. A smooth representation $\pi$ of $\GL_2(F)$ over $\Fpbar$ is called admissible if $\pi^H$ is finite-dimensional for every compact open subgroup $H$ of $\GL_2(F)$. A natural question in the mod $p$ representation theory of $p$-adic groups is whether every smooth irreducible representation is admissible. For $\GL_2(\Qp)$, the answer is affirmative: Berger's theorem guarantees the existence of a central character~\cite{Ber12}, and the classification of Barthel--Livn\'e and Breuil then gives admissibility~\cite{BL94,Bre03}.

For proper finite extensions of $\Qp$, automatic admissibility can fail. For $p>2$, Le constructed non-admissible smooth irreducible mod $p$ representations of $\GL_2(F)$ for unramified extensions $F/\Qp$ of degree at least three~\cite{Le19}. Ghate and Sheth subsequently treated the unramified quadratic extension, also for $p>2$~\cite{GS20}. Ghate, Le, and Sheth extended these results to arbitrary non-archimedean local fields whose residue field is a proper finite extension of $\Fp$, assuming $p>3$~\cite[Theorem~1.1]{GLS23}. In fact, their representations can be chosen to be absolutely irreducible and defined over the residue field of $F$.

The present paper addresses the totally ramified case, in which the residue field is $\Fp$. This case is not covered by the results cited above. For $p\geq 5$, we construct non-admissible smooth irreducible mod $p$ representations of $\GL_2(F)$ when $F/\Qp$ is nontrivial and totally ramified. Thus automatic admissibility can fail even without enlarging the residue field.

Our main result is the following.

\begin{theorem}[Main theorem]\label{thm:main}
Let $p\geq 5$ be a prime, and let $F/\Qp$ be a finite totally ramified extension with $[F:\Qp]>1$. Then $\GL_2(F)$ admits a non-admissible smooth irreducible representation over $\Fpbar$.
\end{theorem}

Our strategy follows the diagram-theoretic approach of the earlier constructions, particularly that of Ghate, Le, and Sheth~\cite{GLS23}. We couple countably many copies of local representations and use injective envelopes to realize the resulting diagram in a smooth representation. The new local input, constructed in Section~\ref{sec:local-constructions}, consists of representations of $\GL_2(\mathcal O_F)$ with nontrivial first congruence subgroup action. Their construction uses nontrivial ramification and replaces the spliced finite-group modules of~\cite{GLS23}. We then deform the operator representing $\Pi$ to couple adjacent summands. Section~\ref{sec:diagram-realization} realizes the modified diagram while controlling the $K$-socle. Finally, Section~\ref{sec:irreducibility} adapts the weighted-shift and minimal-support argument of~\cite[Proposition~3.1]{GLS23} to prove irreducibility of a cyclic subrepresentation. The same shifts produce infinitely many linearly independent vectors fixed by a compact open subgroup, proving non-admissibility.

\paragraph{Related work.}
In independent recent work, Le and Singh~\cite{LS26} also construct non-admissible smooth irreducible mod $p$ representations of $\GL_2(F)$ in the totally ramified case considered here. Both approaches use diagrams with nontrivial $K_1$-action, but the constructions are otherwise distinct. Their work also treats the equal-characteristic case and obtains results for $\GL_n(F)$, for all $n>1$, using parabolic induction.

\paragraph{AI Statement.}
The mathematical construction in Section~\ref{sec:local-constructions} was developed primarily by ChatGPT (OpenAI) through iterative exchanges with the author. This included the local $K$-representation, the compatible $\Pi$-action, and the deformation of the countable direct sum. ChatGPT was also used to draft and revise proofs, identify relevant references, and edit the language and organization of the manuscript. The author is responsible for the mathematical claims and the final manuscript.

\paragraph{Acknowledgements.}
The author is deeply grateful to his advisor, Professor Yongquan Hu, for his guidance, encouragement, and many helpful discussions. The author also thanks Daniel Le and Mansimar Singh for sharing their manuscript~\cite{LS26} and for their interest in this work.

\section{Notation and conventions}\label{sec:notation}

From now on, let $p\geq 5$ be a prime and let $F/\Qp$ be a finite totally ramified extension of degree $e>1$. Write $\mathcal{O}_F$ for its ring of integers, fix a uniformizer $\varpi$, and put $\mathfrak{p}_F=\varpi\mathcal{O}_F$. The residue field $k=\mathcal{O}_F/\mathfrak{p}_F$ is identified with $\Fp$. For $\lambda\in\Fp$, write $[\lambda]\in\mathcal{O}_F$ for its Teichm\"uller lift. Unless otherwise specified, all representations are smooth representations on $\Fpbar$-vector spaces, and all tensor products are over $\Fpbar$.

Set
\[
G=\GL_2(F),\qquad K=\GL_2(\mathcal{O}_F),\qquad
Z=\{z\Id_2:z\in F^\times\}.
\]
Let $\Gamma=\GL_2(\Fp)$ and let $\mathrm{red}:K\twoheadrightarrow\Gamma$ be reduction modulo $\varpi$. Its kernel is the first principal congruence subgroup
\[
K_1=1+\varpi\Mat_2(\mathcal{O}_F).
\]
Let $B$ and $U$ be the upper triangular and upper triangular unipotent subgroups of $\Gamma$, respectively. The corresponding Iwahori subgroup $I$ and the pro-$p$ Iwahori subgroup $I_1$ are
\begin{align*}
I=\mathrm{red}^{-1}(B)
&=\begin{pmatrix}
\mathcal{O}_F^\times & \mathcal{O}_F\\
\mathfrak{p}_F & \mathcal{O}_F^\times
\end{pmatrix},\\
I_1=\mathrm{red}^{-1}(U)
&=\begin{pmatrix}
1+\mathfrak{p}_F & \mathcal{O}_F\\
\mathfrak{p}_F & 1+\mathfrak{p}_F
\end{pmatrix}.
\end{align*}
We have $K/K_1\cong\Gamma$, $I/K_1\cong B$, and $I_1/K_1\cong U$. The subgroup
\[
H=\{\operatorname{diag}([\lambda],[\mu]):\lambda,\mu\in\Fp^\times\}
\]
identifies with $I/I_1$, and $I=H\ltimes I_1$. We also set
\[
s=\begin{pmatrix}0&1\\1&0\end{pmatrix},\qquad
\Pi=\begin{pmatrix}0&1\\\varpi&0\end{pmatrix}.
\]
Then $\Pi$ normalizes both $I$ and $I_1$, and $\Pi^2=\varpi\Id_2$.
Set $N=\langle I,\Pi\rangle$.

A Serre weight is an irreducible $\Fpbar$-representation of $\Gamma$, also viewed as a representation of $K$ by inflation. Every Serre weight is isomorphic to a unique representation of the form
\[
\sigma_{r,m}=\Sym^r\Fpbar^2\otimes\det\nolimits^{m},
\qquad 0\leq r\leq p-1,\quad 0\leq m\leq p-2,
\]
where $\Fpbar^2$ carries the standard action of $\Gamma$. We also write $\sigma_{r,m}=r\otimes\det^m$, interpreting determinant exponents modulo $p-1$; see~\cite[\S1]{BL94}.

The space $\sigma_{r,m}^{I_1}$ is one-dimensional. The action of $I$ on this space factors through $H$ and is given by the character $\chi_{r,m}$:
\[
\chi_{r,m}\bigl(\operatorname{diag}([\lambda],[\mu])\bigr)
=\lambda^{r+m}\mu^m.
\]
For a character $\chi$ of $H$, put $\chi^s(h)=\chi(shs^{-1})$ and also regard $\chi$ as a character of $I$ via $I/I_1\cong H$. For a smooth $K$-representation $V$, write $\soc_K V$ for the sum of its irreducible $K$-subrepresentations, and write $(V^{I_1})^\chi$ for the $\chi$-eigenspace in $V^{I_1}$.

We will use the following generators modulo the pro-$p$ part of the center. Put
\[
u=\begin{pmatrix}1&1\\0&1\end{pmatrix},\qquad
\ell=\begin{pmatrix}1&0\\\varpi&1\end{pmatrix},\qquad
Z_1=(1+\mathfrak p_F)\Id_2.
\]

\begin{lemma}\label{lem:I1-generators}
Let $S=I_1\cap\mathrm{SL}_2(F)$. Then
\[
I_1=Z_1\times S,\qquad S=\overline{\langle u,\ell\rangle},
\]
where the bar denotes closure in $I_1$.
\end{lemma}

\begin{proof}
Since $p\ne2$, taking square roots of determinants in $1+\mathfrak p_F$ gives $I_1=Z_1\times S$. By~\cite[Corollary~5.3]{Koz18}, applied to $\mathrm{SL}_2(F)$, the images of $u$ and $\ell$ form a basis of the Frattini quotient of $S$, since $\mathcal O_F/\mathfrak p_F=\Fp$. The pro-$p$ Burnside basis theorem therefore gives $S=\overline{\langle u,\ell\rangle}$.
\end{proof}

\section{Local constructions and compatible actions}\label{sec:local-constructions}

Fix integers $r$ and $m$ with $2\leq r\leq p-3$ and $0\leq m\leq p-2$. We retain these parameters throughout the construction, until the proof of Theorem~\ref{thm:main}.

\subsection{The local \texorpdfstring{$K$}{K}-representation}

Put $\chi=\chi_{r,m}$ and consider
\[
P_{r,m}=\Ind_I^K\chi^s\cong\Ind_B^\Gamma\chi^s.
\]
$P_{r,m}$ can be viewed as the space of locally constant functions satisfying $f(ik)=\chi^s(i)f(k)$, with right translation $(g\cdot f)(k)=f(kg)$. This representation has dimension $p+1$, is trivial on $K_1$, and fits into the nonsplit exact sequence of~\cite[Theorem~2.4]{BP12} 
\begin{equation}\label{eq:principal-series-extension}
0\longrightarrow\sigma_{r,m}\longrightarrow P_{r,m}
\longrightarrow\sigma_{p-1-r,m+r}\longrightarrow0,
\qquad \soc_K P_{r,m}=\sigma_{r,m}.
\end{equation}

Let $t,\phi\in P_{r,m}^{I_1}$ be supported on $I$ and $IsI_1$, respectively, and normalized by $t(1)=\phi(s)=1$. Then
\begin{equation}\label{eq:principal-series-invariants}
\begin{aligned}
P_{r,m}^{I_1}&=\Fpbar t\oplus\Fpbar\phi,
\qquad \sigma_{r,m}^{I_1}=\Fpbar\phi,\\
h\cdot t&=\chi^s(h)t,\qquad
h\cdot\phi=\chi(h)\phi\quad(h\in H),
\end{aligned}
\end{equation}
and the image of $t$ spans the $I_1$-fixed line of the quotient. We will use the identity
\begin{equation}\label{eq:principal-series-t-to-phi}
\phi=\sum_{\lambda\in\Fp}
\begin{pmatrix}[\lambda]&1\\1&0\end{pmatrix}\cdot t,
\end{equation}
as well as $P_{r,m}=\Fpbar[K]\cdot t$ and $\sigma_{r,m}=\Fpbar[K]\cdot\phi$.

We first construct a $K$-representation $M_{r,m}$, whose underlying $\Fpbar$-vector space is
\begin{equation}\label{eq:module-M-vector-space}
M_{r,m}=P_{r,m}\oplus\sigma_{r+2,m-1}\oplus\sigma_{r,m}.
\end{equation}
The direct sum here is taken as vector spaces. Put $K_2=1+\varpi^2\Mat_2(\mathcal{O}_F)$. Since $e>1$, sending $\varpi$ to $\varepsilon$ identifies $\mathcal{O}_F/(\varpi^2)$ with $\Fp[\varepsilon]/(\varepsilon^2)$, so
\[
K/K_2\cong\GL_2\bigl(\Fp[\varepsilon]/(\varepsilon^2)\bigr)
\cong\Mat_2(\Fp)\rtimes\Gamma,
\]
where the additive group $\Mat_2(\Fp)$ is identified with $1+\varepsilon\Mat_2(\Fp)$ and $\Gamma$ acts on it by conjugation.

We first define the action of $\Mat_2(\Fp)$. We realize $\sigma_{n,\ell}$ as homogeneous polynomials of degree $n$ in $x,y$, with action
\[
\begin{pmatrix}a&b\\c&d\end{pmatrix}\cdot f(x,y)
=(ad-bc)^\ell f(ax+cy,bx+dy).
\]
Fix the embedding
\[
\iota:\sigma_{r,m}\hookrightarrow P_{r,m},\qquad \iota(x^r)=\phi.
\]
For $A=\begin{pmatrix}a&b\\c&d\end{pmatrix}\in\Mat_2(\Fp)$, define
\begin{align}
(1+\varepsilon A)\cdot(v,0,0)&=(v,0,0),\label{eq:module-M-action-P}\\
(1+\varepsilon A)\cdot(0,q,0)
&=\left(\iota\!\left(bq_{yy}+(a-d)q_{xy}-cq_{xx}\right),q,0\right),\label{eq:module-M-action-q}\\
(1+\varepsilon A)\cdot(0,0,w)
&=\left(\iota\!\left(\frac{a-d}{2}(xw_x-yw_y)+bxw_y+cyw_x\right),0,w\right),\label{eq:module-M-action-w}
\end{align}
where $v\in P_{r,m}$, $q\in\sigma_{r+2,m-1}$, $w\in\sigma_{r,m}$, and subscripts denote partial derivatives.

Next, let $\Gamma$ act componentwise by the original actions:
\[
\gamma\cdot(v,q,w)=(\gamma v,\gamma q,\gamma w)
\qquad(\gamma\in\Gamma).
\]
\begin{lemma}\label{lem:module-M-well-defined}
The actions above define a smooth $K$-representation on $M_{r,m}$, trivial on $K_2Z_1$, whose restriction to the subspace $P_{r,m}$ is the original principal series action.
\end{lemma}

\begin{proof}
For $A\in\Mat_2(\Fp)$, let $Q_A(q)$ and $W_A(w)$ denote the polynomial expressions inside $\iota$ in~\eqref{eq:module-M-action-q} and~\eqref{eq:module-M-action-w}, respectively. Put
\[
N_A(v,q,w)=\bigl(\iota(Q_A(q)+W_A(w)),0,0\bigr).
\]
The map $A\mapsto N_A$ is additive, and $N_AN_B=0$ for all $A,B$, since $N_A$ vanishes on $P_{r,m}$. Consequently,
\[
(\Id+N_A)(\Id+N_B)=\Id+N_{A+B},
\]
which proves that these operators define an action of the additive group $\Mat_2(\Fp)$.

To verify compatibility with $\Gamma$, write $A_0=A-\tfrac12\operatorname{tr}(A)\Id_2$ and put
\[
J_0=\begin{pmatrix}0&1\\-1&0\end{pmatrix},\qquad
C_A=-J_0A_0.
\]
The matrix $C_A$ is symmetric, and
\[
Q_A(q)=(\partial_x,\partial_y)C_A
\begin{pmatrix}\partial_x\\\partial_y\end{pmatrix}q.
\]
For $\gamma\in\Gamma$, the identity $\gamma^{\mathrm t}J_0\gamma=\det(\gamma)J_0$ gives
\[
\gamma^{\mathrm t}C_{\gamma A\gamma^{-1}}\gamma=\det(\gamma)C_A.
\]
The chain rule therefore yields
\[
Q_{\gamma A\gamma^{-1}}(\gamma\cdot q)=\gamma\cdot Q_A(q),
\]
where the determinant twists on the source and target are $m-1$ and $m$, respectively. Likewise, $W_A$ is the differential operator induced by $A_0$ on degree-$r$ polynomials, so the chain rule gives
\[
W_{\gamma A\gamma^{-1}}(\gamma\cdot w)=\gamma\cdot W_A(w).
\]
Since $\iota$ is $\Gamma$-equivariant, these identities imply $\gamma N_A\gamma^{-1}=N_{\gamma A\gamma^{-1}}$. Thus the two actions combine to define a representation of $K/K_2\cong\Mat_2(\Fp)\rtimes\Gamma$. Inflation gives the required smooth $K$-representation, with the original action on $P_{r,m}$. For scalar $A$, both $Q_A$ and $W_A$ vanish, so $Z_1$ also acts trivially.
\end{proof}

\begin{lemma}[Socle and extensions]\label{lem:module-M-extensions}
We have
\[
M_{r,m}^{K_1}=P_{r,m},\qquad
\soc_K M_{r,m}=\iota(\sigma_{r,m}).
\]
The $K$-subrepresentations
\[
E_{r,m}=\iota(\sigma_{r,m})\oplus0\oplus\sigma_{r,m},
\qquad
E'_{r,m}=\iota(\sigma_{r,m})\oplus\sigma_{r+2,m-1}\oplus0
\]
define nonzero extension classes
\[
0\ne[E_{r,m}]\in\operatorname{Ext}^1_K(\sigma_{r,m},\sigma_{r,m}),
\qquad
0\ne[E'_{r,m}]\in\operatorname{Ext}^1_K(\sigma_{r+2,m-1},\sigma_{r,m}).
\]
\end{lemma}

\begin{proof}
Let $Q_A$ and $W_A$ denote the differential expressions used in the preceding proof. Consider the $\Gamma$-equivariant map
\[
\begin{aligned}
\sigma_{r+2,m-1}\oplus\sigma_{r,m}
&\longrightarrow\Hom_{\Fp}\bigl(\Mat_2(\Fp),\sigma_{r,m}\bigr),\\
(q,w)&\longmapsto\bigl[A\longmapsto Q_A(q)+W_A(w)\bigr],
\end{aligned}
\]
where $(\gamma\cdot f)(A)=\gamma\cdot f(\gamma^{-1}A\gamma)$ on the target. Its restrictions to the two source summands are nonzero: for the matrix unit $E_{21}$ and $D=\operatorname{diag}(1,-1)$,
\[
Q_{E_{21}}(x^{r+2})=-(r+2)(r+1)x^r,\qquad
W_D(x^r)=r x^r.
\]
These coefficients are nonzero because $2\leq r\leq p-3$. Since the source is the direct sum of two nonisomorphic irreducible $\Gamma$-representations, the map is injective. The defining action of $K_1$ therefore gives $M_{r,m}^{K_1}=P_{r,m}$.

The normal finite $p$-subgroup $K_1/K_2$ acts trivially on every irreducible $K$-subrepresentation of $M_{r,m}$. It follows that
\[
\soc_K M_{r,m}=\soc_K P_{r,m}=\iota(\sigma_{r,m}).
\]
Both $E_{r,m}$ and $E'_{r,m}$ have nontrivial $K_1$-action, whereas a direct sum of their two Serre-weight constituents has trivial $K_1$-action. Their defining extensions are therefore nonsplit.
\end{proof}

The quotient of $M_{r,m}$ by its socle is
\[
M_{r,m}/\iota(\sigma_{r,m})\cong
\sigma_{p-1-r,m+r}\oplus\sigma_{r+2,m-1}\oplus\sigma_{r,m}.
\]
It is semisimple, so $M_{r,m}$ has Loewy length two, with diagram
\[
\begin{tikzcd}[row sep=small, column sep=small, cells={nodes={font=\small}}]
\sigma_{p-1-r,m+r} \arrow[dr, no head] & \sigma_{r+2,m-1} \arrow[d, no head] & \sigma_{r,m} \arrow[dl, no head] \\
& \sigma_{r,m} &
\end{tikzcd}
\]
The upper row is the head, the lower row is the socle, and the lines indicate nontrivial extensions.

\subsection{An \texorpdfstring{$I$}{I}-subrepresentation and the \texorpdfstring{$\Pi$}{Pi}-action}

We next construct the $I$-subrepresentations used to define the $N$-representation $D_1$. Identify $\sigma_{r,m}$ with its image under $\iota$ in the first component of $M_{r,m}$, and write $t$ for $(t,0,0)$. Define
\begin{equation}\label{eq:D0-vectors}
\begin{aligned}
v_0&=(x^r,0,0),&
v_1&=(x^{r-1}y,0,0),&
w_0&=\left(0,0,\frac{x^r}{r}\right),\\
q_0&=\left(0,-\frac{x^{r+2}}{(r+1)(r+2)},0\right),&
q_1&=\left(0,-\frac{x^{r+1}y}{(r+1)(r+2)},0\right).
\end{aligned}
\end{equation}

Let
\begin{equation}\label{eq:module-L-definition}
L_{r,m}=\Fpbar v_0\oplus\Fpbar v_1\oplus\Fpbar w_0
\oplus\Fpbar q_0\oplus\Fpbar q_1\oplus\Fpbar t
\subset M_{r,m}.
\end{equation}

Put $h_0=\operatorname{diag}(1+\varpi,1-\varpi)$ and write
\[
\delta_u=u-1,
\qquad
\delta_{\ell}=\ell-1,
\qquad
\delta_{h_0}=h_0-1
\]
for operators in $\Fpbar[K]$. The formulas~\eqref{eq:module-M-action-P}--\eqref{eq:module-M-action-w} and~\eqref{eq:D0-vectors} give
\begin{equation}\label{eq:D0-operator-actions}
\renewcommand{\arraystretch}{1.5}
\begin{array}{c|ccc}
 & \delta_u & \delta_{\ell} & \delta_{h_0} \\
\hline
v_0 & 0 & 0 & 0 \\
v_1 & v_0 & 0 & 0 \\
q_0 & 0 & v_0 & 0 \\
q_1 & q_0 & \dfrac{r}{r+2}v_1 & -\dfrac{2}{r+2}v_0 \\
w_0 & 0 & v_1 & v_0 \\
t & 0 & 0 & 0
\end{array}
\end{equation}
Indeed, $\delta_u$ acts on each polynomial component by $f(x,y)\mapsto f(x,x+y)-f(x,y)$; for $\delta_{\ell}$ and $\delta_{h_0}$, substitute $A=\begin{pmatrix}0&0\\1&0\end{pmatrix}$ and $A=\operatorname{diag}(1,-1)$, respectively. The last row follows from $t\in P_{r,m}^{I_1}$.

The table shows that $u$ and $\ell$ preserve $L_{r,m}$. Since $Z_1$ acts trivially on $M_{r,m}$, Lemma~\ref{lem:I1-generators} and smoothness imply that $L_{r,m}$ is $I_1$-stable. Each of its basis vectors is an $H$-eigenvector, so $L_{r,m}$ is an $I$-subrepresentation. The same table also gives
\[
\delta_{h_0}=\delta_u\delta_{\ell}-\delta_{\ell}\delta_u
\qquad\text{on }L_{r,m}.
\]

Put $\eta(\operatorname{diag}([a],[b]))=ab^{-1}$ for $a,b\in\Fp^\times$, and retain $\chi=\chi_{r,m}$. In the following Loewy diagram of $L_{r,m}$ as an $I$-representation, each label $z:\theta$ means that $h\cdot z=\theta(h)z$ for $h\in H$:
\[
\begin{tikzcd}[row sep=small, column sep=small, cells={nodes={font=\small}}]
q_1:\chi \arrow[d, no head] \arrow[drr, no head] & & w_0:\chi \arrow[d, no head] & \\
q_0:\chi\eta \arrow[dr, no head] & & v_1:\chi\eta^{-1} \arrow[dl, no head] & \\
& v_0:\chi & \oplus & t:\chi^s
\end{tikzcd}
\]

Define
\begin{align}
D_0&=M_{r,m}\oplus M_{p-1-r,m+r},\label{eq:D0-definition}\\
D_1&=L_{r,m}\oplus L_{p-1-r,m+r}.\label{eq:D1-definition}
\end{align}
The actions on the summands make $D_0$ a $K$-representation and $D_1$ an $I$-representation. The inclusions~\eqref{eq:module-L-definition} induce a natural $I$-equivariant embedding $D_1\hookrightarrow D_0$. 

Put $r'=p-1-r$. We use $v_0,v_1,w_0,q_0,q_1,t$ for the basis of $L_{r,m}$ in~\eqref{eq:module-L-definition}, and primes for the corresponding basis of $L_{r',m+r}$ obtained by replacing $(r,m)$ with $(r',m+r)$.

Define the action of $\Pi$ on the basis of $D_1$ by
\begin{equation}\label{eq:D1-Pi-exchanges}
\begin{aligned}
\Pi v_0&=v'_0,&\qquad \Pi v'_0&=v_0,\\
\Pi v_1&=q'_0,&\qquad \Pi q'_0&=v_1,\\
\Pi q_0&=v'_1,&\qquad \Pi v'_1&=q_0,\\
\Pi t&=t',&\qquad \Pi t'&=t,
\end{aligned}
\end{equation}
and
\begin{equation}\label{eq:D1-Pi-action}
\renewcommand{\arraystretch}{1.2}
\begin{array}{r@{\;}c@{\;}c@{\,}l@{\;}c@{\;}c@{\,}l}
\Pi q_1 & = & \dfrac{r}{r+2} & q'_1 & - & \dfrac{2}{(r+2)(r-1)} & w'_0,\\[6pt]
\Pi w_0 & = & & q'_1 & - & \dfrac{r+1}{r-1} & w'_0,\\[10pt]
\Pi q'_1 & = & \dfrac{r+1}{r-1} & q_1 & - & \dfrac{2}{(r+2)(r-1)} & w_0,\\[6pt]
\Pi w'_0 & = & & q_1 & - & \dfrac{r}{r+2} & w_0.
\end{array}
\end{equation}
Extend these formulas $\Fpbar$-linearly to $D_1$.

\begin{lemma}\label{lem:D1-Pi-compatible}
The prescribed operator $\Pi$ extends the $I$-action on $D_1$ to a smooth action of $\langle I,\Pi\rangle$, with $\varpi\Id_2$ acting trivially.
\end{lemma}

\begin{proof}
Let $P_0$ denote the linear operator prescribed by~\eqref{eq:D1-Pi-exchanges}--\eqref{eq:D1-Pi-action}. These formulas give $P_0^2=\Id$. Using the action table~\eqref{eq:D0-operator-actions}, we obtain
\begin{equation}\label{eq:D1-Pi-commutation}
P_0\delta_u=\delta_{\ell}P_0,\qquad
P_0\delta_{\ell}=\delta_uP_0 \qquad\text{on }D_1.
\end{equation}

Conjugation by the group element $\Pi$ interchanges $u$ and $\ell$ and fixes $Z_1$ pointwise. Since $Z_1$ acts trivially on $D_1$, Lemma~\ref{lem:I1-generators} and smoothness imply
\[
P_0(g\cdot z)=(\Pi g\Pi^{-1})\cdot P_0z
\qquad(g\in I_1,\ z\in D_1).
\]
The operator $P_0$ also exchanges each $H$-character with its conjugate by $s$, in accordance with $\Pi h\Pi^{-1}=shs^{-1}$ for $h\in H$. Since $I=H\ltimes I_1$, the same compatibility holds for every $g\in I$. Thus assigning $P_0$ to $\Pi$ extends the $I$-action to $N$, with $\varpi\Id_2=\Pi^2$ acting trivially. 
\end{proof}
\begin{remark}
The coefficients in~\eqref{eq:D1-Pi-action} are chosen so that $P_0\delta_u=\delta_{\ell}P_0$ and $P_0\delta_{\ell}=\delta_uP_0$ on $L_{r,m}$, and $P_0^2=\Id$. The same intertwining relations then hold on $L_{r',m+r}$. Together with the compatibility of the $H$-characters, these relations give compatibility with the $I$-action.
\end{remark}

\subsection{Deformation of the countable direct sum}

For each $i\in\mathbb{Z}$, let $D_0(i)$ be a copy of the $K$-representation $D_0$ defined in~\eqref{eq:D0-definition}, and let $D_1(i)$ be the corresponding copy of the $I$-representation $D_1$ in~\eqref{eq:D1-definition}. Define
\begin{align}
D_0(\infty)&=\bigoplus_{i\in\mathbb{Z}}D_0(i),\label{eq:D0-infinity}\\
D_1(\infty)&=\bigoplus_{i\in\mathbb{Z}}D_1(i),\label{eq:D1-infinity}
\end{align}
with $K$ acting componentwise on $D_0(\infty)$ and $I$ acting componentwise on $D_1(\infty)$. For $u\in D_j$, write $u_i$ for its image in $D_j(i)$, where $j=0,1$.

The inclusions $D_1(i)\hookrightarrow D_0(i)$ induce a natural $I$-equivariant embedding
\[
D_1(\infty)\hookrightarrow D_0(\infty),
\]
which we use to regard $D_1(\infty)$ as an $I$-subrepresentation of $D_0(\infty)$.

Denote by $P$ the componentwise action of $\Pi$ on $D_1(\infty)$, defined by
\begin{equation}\label{eq:D1-infinity-Pi-action}
P u_i=(\Pi u)_i
\qquad(u\in D_1,\ i\in\mathbb{Z}),
\end{equation}
where the action on the right is defined in~\eqref{eq:D1-Pi-exchanges}--\eqref{eq:D1-Pi-action}. Together with the componentwise $I$-action, $P$ defines an $N$-action on $D_1(\infty)$, with $P^2=\Id$.

We modify $P$ by adding terms in adjacent summands. The coefficients are chosen so that the modified operator remains an involution and is compatible with the $I$-action. Fix a sequence $\lambda=(\lambda_i)_{i\in\mathbb{Z}}$ with $\lambda_i\in\Fpbar$, and define an $\Fpbar$-linear operator $\Delta$ on $D_1(\infty)$ by
\begin{equation}\label{eq:D1-infinity-Delta-action}
\begin{aligned}
\Delta((q_1)_i)&=\lambda_i t_{i+1},\\[4pt]
\Delta((w_0)_i)&=\lambda_i t_{i-1},\\[6pt]
\Delta((q'_1)_i)&=-\frac{r+1}{r-1}\lambda_i t'_{i+1}
+\frac{2}{(r+2)(r-1)}\lambda_i t'_{i-1},\\[6pt]
\Delta((w'_0)_i)&=-\lambda_i t'_{i+1}
+\frac{r}{r+2}\lambda_i t'_{i-1}
\qquad(i\in\mathbb{Z}),
\end{aligned}
\end{equation}
and let $\Delta$ vanish on all the remaining basis vectors $(v_0)_i,(v_1)_i,(q_0)_i,t_i,(v'_0)_i,(v'_1)_i,(q'_0)_i,t'_i$. These formulas define an endomorphism of the direct sum $D_1(\infty)$.

Put
\begin{equation}\label{eq:D1-infinity-modified-Pi}
\widetilde P=P+\Delta\in\End_{\Fpbar}\bigl(D_1(\infty)\bigr).
\end{equation}

\begin{proposition}\label{prop:D1-infinity-modified-action}
For every sequence $\lambda$, the operator $\widetilde P$ satisfies
\[
\widetilde P^2=\Id,\qquad
\widetilde P(g\cdot u)=(\Pi g\Pi^{-1})\cdot\widetilde P u
\qquad(g\in I,\ u\in D_1(\infty)).
\]
Consequently, keeping the $I$-action unchanged and letting the group element $\Pi$ act by $\widetilde P$ defines a smooth representation of $\langle I,\Pi\rangle$ on $D_1(\infty)$, with $\varpi\Id_2$ acting trivially.
\end{proposition}

\begin{proof}
The image of $\Delta$ is contained in
\[
T=\bigoplus_{i\in\mathbb{Z}}\bigl(\Fpbar t_i\oplus\Fpbar t'_i\bigr),
\]
and $\Delta$ vanishes on $T$. Thus $\Delta^2=0$. Substituting~\eqref{eq:D1-Pi-exchanges}--\eqref{eq:D1-Pi-action} into~\eqref{eq:D1-infinity-Delta-action} gives
\[
\Delta P((q_1)_i)=-\lambda_i t'_{i+1},\qquad
\Delta P((w_0)_i)=-\lambda_i t'_{i-1}.
\]
These are $-P\Delta((q_1)_i)$ and $-P\Delta((w_0)_i)$, respectively. For the primed vectors, the last two formulas in~\eqref{eq:D1-infinity-Delta-action} give the same identity directly. Both compositions vanish on the remaining basis vectors. Hence
\begin{equation}\label{eq:D1-infinity-Delta-identities}
\Delta^2=0,\qquad P\Delta+\Delta P=0,
\end{equation}
and therefore $\widetilde P^2=P^2+P\Delta+\Delta P+\Delta^2=\Id$.

To check compatibility with $I_1$, the action table~\eqref{eq:D0-operator-actions} and the definition of $\Delta$ give
\[
\Delta\delta_u=\Delta\delta_{\ell}
=\delta_u\Delta=\delta_{\ell}\Delta=0,
\]
where $\delta_u$ and $\delta_{\ell}$ act componentwise on $D_1(\infty)$. Thus
\[
\Delta(u\cdot z)=\ell\cdot\Delta z,\qquad
\Delta(\ell\cdot z)=u\cdot\Delta z.
\]
Since $Z_1$ acts trivially and $\Pi$ interchanges $u$ and $\ell$ by conjugation, Lemma~\ref{lem:I1-generators} and smoothness imply
\[
\Delta(g\cdot z)=(\Pi g\Pi^{-1})\cdot\Delta z
\qquad(g\in I_1,\ z\in D_1(\infty)).
\]

For compatibility with $H$, write $\chi=\chi_{r,m}$. The vectors $q_1,w_0$ have $H$-character $\chi$, whereas $q'_1,w'_0$ have character $\chi^s$. The vectors $t$ and $t'$ have characters $\chi^s$ and $\chi$, respectively. Thus every nonzero image in~\eqref{eq:D1-infinity-Delta-action} has the character obtained from that of its source by conjugation by $s$. Since $\Pi h\Pi^{-1}=shs^{-1}$ for $h\in H$, this proves the same compatibility for $H$. Using $I=H\ltimes I_1$, we obtain it for all $g\in I$. Adding the corresponding identity for $P$ proves the assertion for $\widetilde P$.
\end{proof}

\Needspace{9\baselineskip}
\section{From diagrams to smooth representations}\label{sec:diagram-realization}

\subsection{Diagrams and injective envelopes}

\begin{definition}[{\cite[Definition~9.7]{BP12}}]\label{def:diagrams}
A diagram is a triple $(V_0,V_1,\jmath)$, where $V_0$ is a smooth $KZ$-representation, $V_1$ is a smooth $N$-representation, and $\jmath:V_1\to V_0$ is an $IZ$-equivariant map. A morphism of diagrams is a pair of equivariant maps on the two components that commutes with the connecting maps. It is an embedding if both component maps are injective.
\end{definition}

Extend the $K$-actions on $D_0$ and $D_0(\infty)$ to $KZ$ by letting $\varpi\Id_2$ act trivially. On $D_1$, retain the action of $\langle I,\Pi\rangle$ constructed in Lemma~\ref{lem:D1-Pi-compatible}. The natural inclusion $\operatorname{incl}:D_1\hookrightarrow D_0$ is then $IZ$-equivariant: it is $I$-equivariant, and $\varpi\Id_2$ acts trivially on both sides. Define the diagram
\begin{equation}\label{eq:diagram-D}
D=(D_0,D_1,\operatorname{incl}).
\end{equation}

Likewise, define
\begin{equation}\label{eq:diagram-D-infinity}
D(\infty)=\bigl(D_0(\infty),D_1(\infty),\operatorname{incl}\bigr),
\end{equation}
where the group element $\Pi$ acts on $D_1(\infty)$ by the original componentwise operator $P$. Thus $D(\infty)$ is the direct sum, indexed by $\mathbb{Z}$, of copies of the diagram $D$.

Finally, for the fixed sequence $\lambda$, define the modified diagram
\begin{equation}\label{eq:diagram-modified-D-infinity}
\widetilde D_{\lambda}(\infty)=\bigl(D_0(\infty),D_1(\infty),\operatorname{incl}\bigr),
\end{equation}
where $\Pi$ acts on $D_1(\infty)$ by $\widetilde P=P+\Delta$. This is a well-defined diagram by Proposition~\ref{prop:D1-infinity-modified-action}; it differs from $D(\infty)$ only in the operator representing $\Pi$.

For a smooth representation $V$ of $K$ or $I$, write $\Inj_K(V)$ or $\Inj_I(V)$ for its injective envelope in the corresponding category of smooth representations over $\Fpbar$.

\begin{proposition}\label{prop:diagram-injective-embedding}
There exists a smooth $G$-representation $\Omega$, with $\varpi\Id_2$ acting trivially, and an embedding of diagrams
\begin{equation}\label{eq:diagram-D-embedding-Omega}
D\hookrightarrow\bigl(\Omega|_{KZ},\Omega|_N,\Id_{\Omega}\bigr),
\qquad \Omega|_K\cong\Inj_K(D_0).
\end{equation}
After identifying $D_0$ and $D_1$ with their images, we obtain an $I$-equivariant projection $e_I$ satisfying
\begin{equation}\label{eq:Omega-I-projection}
e_I\in\End_I(\Omega),\qquad e_I^2=e_I,\qquad
 e_I|_{D_1}=\Id_{D_1},\qquad e_I(\Omega)\cong\Inj_I(D_1),
\end{equation}
where the inclusions $D_0\hookrightarrow\Omega$ and $D_1\hookrightarrow e_I(\Omega)$ realize the respective injective envelopes.
\end{proposition}

\begin{proof}
This follows from~\cite[Theorem~9.8]{BP12}.
\end{proof}

Fix such an embedding and projection, and identify $D_0,D_1$ with their images in $\Omega$. The projection gives a decomposition of $I$-representations
\[
\Omega=e_I(\Omega)\oplus(1-e_I)(\Omega).
\]

For each $i\in\mathbb{Z}$, let $\Omega(i)$ be a copy of the smooth $G$-representation $\Omega$, and put
\begin{equation}\label{eq:Omega-infinity}
\Omega(\infty)=\bigoplus_{i\in\mathbb{Z}}\Omega(i),
\end{equation}
with $G$ acting componentwise. Taking the direct sum of the fixed embedding~\eqref{eq:diagram-D-embedding-Omega} over all $i\in\mathbb{Z}$ gives an embedding of diagrams
\begin{equation}\label{eq:diagram-D-infinity-embedding-Omega}
D(\infty)\hookrightarrow
\bigl(\Omega(\infty)|_{KZ},\Omega(\infty)|_N,\Id_{\Omega(\infty)}\bigr).
\end{equation}
We use this embedding to identify $D_0(\infty)$ and $D_1(\infty)$ with their images in $\Omega(\infty)$. In particular, the action of $\Pi$ on $\Omega(\infty)$ restricts to the original operator $P$ on $D_1(\infty)$.

For each $i\in\mathbb{Z}$, let $e_{I,i}\in\End_I(\Omega(i))$ be the projection corresponding to $e_I$ under the identification $\Omega(i)\cong\Omega$. Define the $I$-subrepresentation
\begin{equation}\label{eq:calE-definition}
\mathcal E=\bigoplus_{i\in\mathbb{Z}}e_{I,i}\bigl(\Omega(i)\bigr)
\subset\Omega(\infty).
\end{equation}
Thus $\mathcal E$ is the direct sum of the copies of $e_I(\Omega)$ inside $\Omega(\infty)$, and $D_1(\infty)\subset\mathcal E$.

\begin{lemma}\label{lem:calE-injective-envelope}
The inclusion $D_1(\infty)\hookrightarrow\mathcal E$ realizes an injective envelope in the category of smooth $I$-representations:
\begin{equation}\label{eq:calE-injective-envelope}
\mathcal E\cong\Inj_I\bigl(D_1(\infty)\bigr).
\end{equation}
\end{lemma}

\begin{proof}
In the category of smooth $I$-representations, arbitrary direct sums of injective objects are injective: this category is locally Noetherian, since every smooth $I$-representation is the union of its finite-dimensional $I$-subrepresentations. Each $e_{I,i}(\Omega(i))$ is injective, so $\mathcal E$ is injective as a smooth $I$-representation.

Moreover, $D_1(\infty)\hookrightarrow\mathcal E$ is an essential extension. Indeed, any nonzero $I$-subrepresentation of $\mathcal E$ contains a nonzero vector supported on finitely many summands. The $I$-subrepresentation generated by this vector is contained in the corresponding finite direct sum of the $e_{I,i}(\Omega(i))$. Since each $D_1(i)\hookrightarrow e_{I,i}(\Omega(i))$ is essential, their finite direct sum is essential as well. Hence the given nonzero subrepresentation meets $D_1(\infty)$ nontrivially. Consequently, the inclusion realizes $\mathcal E$ as an injective envelope.
\end{proof}

\subsection{Extending the deformation}

On the subspace $D_1(\infty)\subset\Omega(\infty)$, define
\begin{equation}\label{eq:D1-infinity-conjugating-automorphism}
\alpha=\Id+\frac12\Delta P.
\end{equation}
Both $\Delta$ and $P$ intertwine the $I$-action with its conjugate by $\Pi$. Since conjugation by $\Pi$ has order two on $I$, their composite $\Delta P$ is $I$-equivariant. Moreover,~\eqref{eq:D1-infinity-Delta-identities} gives
\[
(\Delta P)^2=-\Delta^2P^2=0.
\]
Thus $\alpha$ is an $I$-equivariant automorphism of $D_1(\infty)$, with
\[
\alpha^{-1}=\Id-\frac12\Delta P.
\]

\begin{lemma}\label{lem:conjugating-automorphism-extension}
The automorphism $\alpha$ extends to an $I$-equivariant automorphism $\widehat\alpha$ of $\Omega(\infty)$ which preserves $\mathcal E$ and acts as the identity on $\bigoplus_{i\in\mathbb{Z}}(1-e_{I,i})(\Omega(i))$.
\end{lemma}

\begin{proof}
By the injectivity of $\mathcal E$, the map $\alpha:D_1(\infty)\to\mathcal E$ extends to an $I$-equivariant endomorphism $\alpha_{\mathcal E}$ of $\mathcal E$. This extension is an automorphism. Indeed, its kernel meets $D_1(\infty)$ trivially, so essentiality implies that its kernel is zero. Its image is then isomorphic to the injective object $\mathcal E$ and hence is a direct summand of $\mathcal E$. Since this image contains $\alpha(D_1(\infty))=D_1(\infty)$, essentiality forces its complementary summand to vanish. Fix such an extension $\alpha_{\mathcal E}$.

The projections $e_{I,i}$ give the decomposition of $I$-representations
\[
\Omega(\infty)=\mathcal E\oplus
\bigoplus_{i\in\mathbb{Z}}(1-e_{I,i})\bigl(\Omega(i)\bigr).
\]
Extending $\alpha_{\mathcal E}$ by the identity on the second summand, we obtain an $I$-equivariant automorphism $\widehat\alpha$ of $\Omega(\infty)$. Its restriction to $D_1(\infty)$ is $\alpha$.
\end{proof}

\begin{proposition}\label{prop:modified-diagram-realization}
For every sequence $\lambda$, there is a smooth $G$-representation $\widetilde\Omega(\infty)$ whose underlying vector space and $KZ$-action agree with those of $\Omega(\infty)$. The fixed inclusions embed $\widetilde D_\lambda(\infty)$ into the diagram associated with $\widetilde\Omega(\infty)$. The $G$-subrepresentation $\pi$ generated by $D_0(\infty)$ satisfies
\[
\soc_K\pi=\soc_K D_0(\infty).
\]
\end{proposition}

\begin{proof}
Let $\rho$ denote the original componentwise $G$-action on $\Omega(\infty)$, and define the new operator representing $\Pi$ by
\begin{equation}\label{eq:Omega-infinity-modified-Pi}
\widetilde\rho(\Pi)=\widehat\alpha\,\rho(\Pi)\,\widehat\alpha^{-1}.
\end{equation}
Since $\widehat\alpha$ is $I$-equivariant, this operator satisfies
\[
\widetilde\rho(\Pi)^2=\Id,\qquad
\widetilde\rho(\Pi)\rho(g)=\rho(\Pi g\Pi^{-1})\widetilde\rho(\Pi)
\qquad(g\in I).
\]
It therefore defines, together with the unchanged $I$-action, a new smooth $N$-action on $\Omega(\infty)$, with $\varpi\Id_2$ acting trivially. On $D_1(\infty)$, we have
\begin{align*}
\widetilde\rho(\Pi)|_{D_1(\infty)}
&=\alpha P\alpha^{-1}\\
&=\left(\Id+\frac12\Delta P\right)P
\left(\Id-\frac12\Delta P\right)\\
&=P+\frac12\Delta-\frac12P\Delta P-\frac14\Delta^2P\\
&=P+\Delta=\widetilde P,
\end{align*}
where we used $P\Delta P=-\Delta$ and $\Delta^2=0$.

The original $KZ$-action and the new $N$-action on $\Omega(\infty)$ agree on $IZ$: the $I$-action is unchanged, and $\varpi\Id_2$ acts trivially in both actions. Hence Pa\v{s}k\=unas~\cite[Corollary~5.5.5]{Pa04} gives a new smooth $G$-representation on this vector space with the original $KZ$-action and modified $\Pi$-action $\tilde{\rho}(\Pi)$. Denote it by $\widetilde\Omega(\infty)$. By construction, there is an embedding
\begin{equation}\label{eq:modified-diagram-embedding-Omega}
\widetilde\iota:\widetilde D_{\lambda}(\infty)\hookrightarrow
\bigl(\widetilde\Omega(\infty)|_{KZ},
\widetilde\Omega(\infty)|_N,\Id_{\widetilde\Omega(\infty)}\bigr).
\end{equation}

Define
\begin{equation}\label{eq:pi-definition}
\pi=\Fpbar[G]\cdot D_0(\infty)
\subset\widetilde\Omega(\infty).
\end{equation}
Then $\pi$ is smooth, and the fixed inclusion gives
\[
D_0(\infty)\subset\pi\subset\widetilde\Omega(\infty).
\]

Finally, the inclusion $D_0\hookrightarrow\Omega$ is an injective envelope as a $K$-representation, so $\soc_K\Omega=\soc_K D_0$. Taking direct sums, we obtain
\[
\soc_K\widetilde\Omega(\infty)
=\soc_K\Omega(\infty)
=\soc_K D_0(\infty).
\]
The inclusions above therefore imply
\begin{equation}\label{eq:pi-K-socle}
\soc_K\pi=\soc_K D_0(\infty)
\cong\bigoplus_{i\in\mathbb{Z}}
\bigl(\sigma_{r,m}\oplus\sigma_{p-1-r,m+r}\bigr).
\end{equation}
\end{proof}

We retain these choices of $\widetilde\Omega(\infty)$ and $\pi$ in what follows.

\section{Irreducibility and non-admissibility}\label{sec:irreducibility}

The irreducibility argument in this section is modeled on the proof of~\cite[Proposition~3.1]{GLS23}. In that proof, the spliced-module construction yields weighted shifts in both directions on finitely supported coefficient sequences; a minimal-support argument and the condition on $\lambda_{i-1}\lambda_i$ then isolate a single summand. Here, the operators $\mathcal F$, $J$, and the modified $\Pi$-action will yield the analogous shifts $A$ and $B$. Their composite is diagonal with eigenvalues $\lambda_{i-1}\lambda_i$, allowing the same support-reduction argument. We adapt this method to a cyclic $G$-subrepresentation of $\pi$, rather than applying the diagram-irreducibility statement of~\cite[Proposition~3.1]{GLS23} directly.

Let $(v_0)_0$ denote the copy of $v_0$ in $D_0(0)\subset\pi$. Define the cyclic $G$-subrepresentation
\begin{equation}\label{eq:tilde-pi-definition}
\widetilde\pi=\Fpbar[G]\cdot(v_0)_0\subset\pi.
\end{equation}
We will prove the following irreducibility statement.

\begin{proposition}\label{prop:tilde-pi-irreducible}
Suppose that $\lambda_i\in\Fpbar^\times$ for every $i\in\mathbb{Z}$ and that the sequence $\lambda=(\lambda_i)_{i\in\mathbb{Z}}$ satisfies
\begin{equation}\label{eq:lambda-irreducibility-condition}
\lambda_{i-1}\lambda_i\ne\lambda_{-1}\lambda_0
\qquad\text{for every }i\in\mathbb{Z}\setminus\{0\}.
\end{equation}
Then $\widetilde\pi$ is an irreducible smooth $G$-representation.
\end{proposition}

To prove Proposition~\ref{prop:tilde-pi-irreducible}, we introduce two operators in $\Fpbar[K]$:
\begin{equation}\label{eq:F-J-definition}
\mathcal F=\sum_{a\in\Fp}a^{p-2}
\begin{pmatrix}1&0\\{[a]}&1\end{pmatrix},
\qquad
J=\sum_{a\in\Fp}\begin{pmatrix}[a]&1\\1&0\end{pmatrix}.
\end{equation}
\begin{lemma}[Auxiliary operators]\label{lem:auxiliary-operators}
The operators $\mathcal F$ and $J$ satisfy
\begin{equation}\label{eq:F-J-basis-actions}
\begin{aligned}
\mathcal F v_0&=-r v_1, & \mathcal F q_0&=-(r+2)q_1,\\
Jq_1&=Jw_0=0, & Jt&=v_0.
\end{aligned}
\end{equation}
The corresponding formulas for the primed vectors hold with $(r,m)$ replaced by $(r',m+r)$, where $r'=p-1-r$.
\end{lemma}

\begin{proof}
Since $e>1$, the images modulo $\varpi^2$ of the matrices occurring in $\mathcal F$ and $J$ lie in the constant subgroup $\Gamma\subset\GL_2(\Fp[\varepsilon]/(\varepsilon^2))$. They therefore act componentwise on $M_{r,m}$, with no additional terms from its extension structure.

For $0\leq j\leq p-1$, the finite-field power sums give
\[
\sum_{a\in\Fp}a^{p-2+j}=
\begin{cases}-1,&j=1,\\0,&j\ne1.\end{cases}
\]
The matrices defining $\mathcal F$ have determinant one and send $x$ to $x+ay$ and $y$ to $y$. Thus, for $1\leq d\leq p-1$,
\[
\mathcal F(x^d)=\sum_{j=0}^d\binom dj x^{d-j}y^j
\sum_{a\in\Fp}a^{p-2+j}=-d x^{d-1}y.
\]
Taking $d=r$ and $d=r+2$, and using the normalizations in~\eqref{eq:D0-vectors}, proves the two identities for $\mathcal F$.

On $\sigma_{d,n}$, the matrix $\begin{pmatrix}a&1\\1&0\end{pmatrix}$ sends $f(x,y)$ to $(-1)^n f(ax+y,x)$. Consequently,
\begin{align*}
Jq_1&=\left(0,-\frac{(-1)^{m-1}x}{(r+1)(r+2)}
\sum_{a\in\Fp}(ax+y)^{r+1},0\right)=0,\\
Jw_0&=\left(0,0,\frac{(-1)^m}{r}
\sum_{a\in\Fp}(ax+y)^r\right)=0.
\end{align*}
Indeed, $r+1\leq p-2$, and $\sum_{a\in\Fp}a^j=0$ for every $0\leq j\leq p-2$, with $a^0$ interpreted as the constant polynomial $1$. Finally, $Jt=v_0$ is~\eqref{eq:principal-series-t-to-phi}. Since $r'$ also satisfies $2\leq r'\leq p-3$, the same argument proves the primed identities.
\end{proof}

Let
\[
\mathcal C=\bigoplus_{i\in\mathbb{Z}}\Fpbar\mathbf e_i
\]
be the space of finitely supported coefficient sequences, with standard basis $(\mathbf e_i)_{i\in\mathbb{Z}}$. Define operators $A,B\in\End_{\Fpbar}(\mathcal C)$ by
\begin{equation}\label{eq:coefficient-shifts}
A\mathbf e_i=\lambda_i\mathbf e_{i+1},\qquad
B\mathbf e_i=\lambda_i\mathbf e_{i-1}.
\end{equation}
For a $G$-subrepresentation $W\subset\pi$, put
\begin{equation}\label{eq:coefficient-space-W}
\mathcal C_W=\left\{\sum_i c_i\mathbf e_i\in\mathcal C:
\sum_i c_i(v_0)_i\in W\right\}.
\end{equation}
All sums in this definition are finite.

\begin{lemma}\label{lem:coefficient-shift-stability}
For every $G$-subrepresentation $W\subset\pi$, the subspace $\mathcal C_W$ is stable under both $A$ and $B$.
\end{lemma}

\begin{proof}
Throughout the proof, $\Pi$ denotes the group element acting on $\pi$; its restriction to $D_1(\infty)$ is $\widetilde P$, not $P$. Fix $\sum_i c_i\mathbf e_i\in\mathcal C_W$ and put $u=\sum_i c_i(v_0)_i\in W$. The formulas for $\mathcal F$ and the unchanged exchanges of $v_0,v_1,q_0$ and their primed counterparts give
\begin{align*}
\mathcal F\Pi\mathcal F\Pi u
&=r'(r+2)\sum_i c_i(q_1)_i,\\
\mathcal F\Pi\mathcal F u
&=r(r'+2)\sum_i c_i(q'_1)_i.
\end{align*}
Both coefficients are nonzero in $\Fpbar$, since $2\leq r,r'\leq p-3$. Thus both sums on the right belong to $W$.

By~\eqref{eq:D1-infinity-Delta-action} and~\eqref{eq:F-J-basis-actions},
\begin{align}
J\Pi(q_1)_i
&=\lambda_i(v_0)_{i+1},\label{eq:J-Pi-q1}\\
\Pi J\Pi(q'_1)_i
&=-\frac{r+1}{r-1}\lambda_i(v_0)_{i+1}
+\frac{2}{(r+2)(r-1)}\lambda_i(v_0)_{i-1}.
\label{eq:Pi-J-Pi-q1-prime}
\end{align}
Indeed, $J$ kills $q_1,w_0,q'_1,w'_0$, whereas $Jt=v_0$ and $Jt'=v'_0$; moreover, $\Pi(v'_0)_j=(v_0)_j$ for every $j$. Applying~\eqref{eq:J-Pi-q1} to $\sum_i c_i(q_1)_i\in W$ yields
\[
\sum_i c_i\lambda_i(v_0)_{i+1}\in W,
\]
which proves $A\mathcal C_W\subset\mathcal C_W$. Applying~\eqref{eq:Pi-J-Pi-q1-prime} to $\sum_i c_i(q'_1)_i\in W$ and using this first inclusion then gives
\[
\sum_i c_i\lambda_i(v_0)_{i-1}\in W,
\]
since $2/((r+2)(r-1))\ne0$. Hence $B\mathcal C_W\subset\mathcal C_W$ as well.
\end{proof}

This stability holds for every sequence $\lambda$, without assuming~\eqref{eq:lambda-irreducibility-condition}.

\begin{proof}[Proof of Proposition~\ref{prop:tilde-pi-irreducible}]
We prove the stronger assertion that every nonzero $G$-subrepresentation of $\pi$ contains $\widetilde\pi$. Let $W\subset\pi$ be such a subrepresentation. We first show that $\mathcal C_W\ne0$. Since $K$ is compact and $W$ is smooth, the $K$-subrepresentation generated by any nonzero vector of $W$ is finite-dimensional and contains an irreducible $K$-subrepresentation $\tau$. By~\eqref{eq:pi-K-socle}, the representation $\tau$ is isomorphic to either $\sigma_{r,m}$ or $\sigma_{p-1-r,m+r}$. These two Serre weights are nonisomorphic. If $\tau\cong\sigma_{r,m}$, a nonzero vector of $\tau^{I_1}$ has the form
\[
\sum_i c_i(v_0)_i\in W
\]
for a nonzero finitely supported sequence $(c_i)$. If $\tau\cong\sigma_{p-1-r,m+r}$, such a vector instead has the form $\sum_i c_i(v'_0)_i$; applying $\Pi$ gives $\sum_i c_i(v_0)_i\in W$. In either case, $\mathcal C_W\ne0$.

Choose a nonzero vector
\[
v=\sum_i c_i\mathbf e_i\in\mathcal C_W
\]
for which $|\operatorname{supp}(v)|$ is minimal, where $\operatorname{supp}(v)=\{i\in\mathbb{Z}:c_i\ne0\}$. By Lemma~\ref{lem:coefficient-shift-stability}, the space $\mathcal C_W$ is stable under $A$ and $B$. Since every $\lambda_i$ is nonzero, applying either operator shifts the support by one without changing its cardinality. Choose $j\in\operatorname{supp}(v)$. Since $\lambda_i \neq 0$, replacing $v$ by $B^jv$ if $j>0$, or by $A^{-j}v$ if $j<0$, we may therefore assume that $c_0\ne0$. 

The composite $AB$ is diagonal on $\mathcal C$:
\[
AB\mathbf e_i=\lambda_{i-1}\lambda_i\mathbf e_i.
\]
Consequently,
\begin{equation}\label{eq:minimal-support-AB}
\bigl(AB-\lambda_{-1}\lambda_0\Id\bigr)v
=\sum_i\bigl(\lambda_{i-1}\lambda_i-\lambda_{-1}\lambda_0\bigr)c_i\mathbf e_i
\in\mathcal C_W.
\end{equation}
The coefficient at $i=0$ vanishes, and the support of this vector is contained in $\operatorname{supp}(v)\setminus\{0\}$. By minimality,
\[
\bigl(AB-\lambda_{-1}\lambda_0\Id\bigr)v=0.
\]
By~\eqref{eq:lambda-irreducibility-condition}, each factor $\lambda_{i-1}\lambda_i-\lambda_{-1}\lambda_0$ is nonzero for $i\ne0$, so $c_i=0$ for every $i\ne0$. Thus $v=c_0\mathbf e_0$ with $c_0\ne0$, and $\mathbf e_0\in\mathcal C_W$. By definition, $(v_0)_0\in W$ and hence
\[
\widetilde\pi=\Fpbar[G]\cdot(v_0)_0\subset W.
\]

We have shown that every nonzero $G$-subrepresentation of $\pi$ contains $\widetilde\pi$. In particular, any nonzero $G$-subrepresentation of $\widetilde\pi$ equals $\widetilde\pi$. This proves its irreducibility.
\end{proof}

\Needspace{9\baselineskip}
\begin{corollary}\label{cor:tilde-pi-non-admissible}
Under the hypotheses of Proposition~\ref{prop:tilde-pi-irreducible}, the representation $\widetilde\pi$ is non-admissible.
\end{corollary}

\begin{proof}
Since $(v_0)_0\in\widetilde\pi$, we have $\mathbf e_0\in\mathcal C_{\widetilde\pi}$. By Lemma~\ref{lem:coefficient-shift-stability} and the nonvanishing of every $\lambda_i$, repeated applications of $A$ and $B$ show that $\mathbf e_i\in\mathcal C_{\widetilde\pi}$ for all $i\in\mathbb{Z}$. Thus $(v_0)_i\in\widetilde\pi$ for every $i$. These vectors are linearly independent, since they lie in distinct summands of $D_0(\infty)$, and each is fixed by $K_1$, since $v_0$ belongs to the inflated Serre weight $\sigma_{r,m}$. Hence
\[
\bigoplus_{i\in\mathbb{Z}}\Fpbar(v_0)_i
\subset\widetilde\pi^{K_1},\qquad
\dim_{\Fpbar}\widetilde\pi^{K_1}=\infty.
\]
As $K_1$ is compact and open in $G$, the representation $\widetilde\pi$ is not admissible.
\end{proof}

\begin{proof}[Proof of Theorem~\ref{thm:main}]
For every $p\geq5$, we may take $r=2$ and $m=0$ in the preceding construction. Choose the sequence $\lambda$ by
\[
\lambda_i=1\quad(i\geq0),\qquad
\lambda_{-n}=2^{(-1)^{n+1}}\quad(n\geq1),
\]
where $2$ is viewed as an element of $\Fp^\times$. Every $\lambda_i$ is nonzero, and
\[
\lambda_{i-1}\lambda_i=
\begin{cases}
2,&i=0,\\
1,&i\ne0.
\end{cases}
\]
Thus~\eqref{eq:lambda-irreducibility-condition} holds. Proposition~\ref{prop:tilde-pi-irreducible} and Corollary~\ref{cor:tilde-pi-non-admissible} now show that $\widetilde\pi$ is a non-admissible smooth irreducible representation of $\GL_2(F)$ over $\Fpbar$, as required.
\end{proof}

\begin{remark}\label{rem:tilde-pi-socle-and-pi}
Under the hypotheses of Proposition~\ref{prop:tilde-pi-irreducible}, the representations $\widetilde\pi$ and $\pi$ have the same $K$-socle:
\[
\soc_K\widetilde\pi=\soc_K\pi=\soc_K D_0(\infty)
\cong\bigoplus_{i\in\mathbb Z}
\bigl(\sigma_{r,m}\oplus\sigma_{p-1-r,m+r}\bigr).
\]
Indeed, the proof of Corollary~\ref{cor:tilde-pi-non-admissible} shows that $(v_0)_i\in\widetilde\pi$ for every $i$. Applying $\Pi$ also gives $(v'_0)_i\in\widetilde\pi$. Their $K$-spans exhaust $\soc_K D_0(\infty)$, so the equality follows from~\eqref{eq:pi-K-socle} and $\widetilde\pi\subset\pi$.

Moreover, the proof of Proposition~\ref{prop:tilde-pi-irreducible} shows that every nonzero $G$-subrepresentation of $\pi$ contains $\widetilde\pi$. Thus $\widetilde\pi$ is the unique irreducible $G$-subrepresentation of $\pi$. Equality $\widetilde\pi=\pi$ is equivalent to $D_0(\infty)\subset\widetilde\pi$, which is not established by the preceding argument.
\end{remark}

\end{document}